\documentclass[12pt]{article}
\usepackage{amsmath}
\usepackage{amssymb}
\begin{document}
\newtheorem{theorem}{Theorem}[section]
\newtheorem{claim}[theorem]{Claim}
\newtheorem{tlemma}[theorem]{Test Lemma}
\newtheorem{llemma}[theorem]{Lifting Lemma}
\newtheorem{remark}[theorem]{Remark}
\newtheorem{definition}[theorem]{Definition}
\newtheorem{notation}[theorem]{Notation}
\newtheorem{Sh-ele}[theorem]{Shelah's Elevator}
\newtheorem{mtheorem}[theorem]{Main Theorem}
\newtheorem{steplemma}[theorem]{Step Lemma}
\newtheorem{blackbox}[theorem]{The General Black Box}
\newtheorem{weakdia}[theorem]{The Weak Diamond Principle }
\newtheorem{sblackbox}[theorem]{The Strong Black Box}
\newtheorem{wllemma}[theorem]{Wald-{\L}o{\'s} Lemma}
\newtheorem{pigeonlemma}[theorem]{Pigeon-hole Lemma}
\newtheorem{observation}[theorem]{Observation}
\newtheorem{proposition}[theorem]{Proposition}
\newtheorem{lemma}[theorem]{Lemma}
\newtheorem{corollary}[theorem]{Corollary}
\newtheorem{construction}[theorem]{Construction}
\newtheorem{example}[theorem]{Example}
\newtheorem{recognition}[theorem]{Recognition Lemma}
\renewcommand{\labelenumi}{(\roman{enumi})}
\numberwithin{equation}{section}
\def\Pf{\smallskip\goodbreak{\sl Proof. }}
\def\Fp{\vadjust{}\penalty200 \hfill
\lower.3333ex\hbox{\vbox{\hrule\hbox{\vrule\phantom{\vrule height
6.83333pt depth 1.94444pt width 8.77777pt}\vrule}\hrule}}
\ifmmode\let\next\relax\else\let\next\par\fi \next}
\def\qa{(A_\a, G_\a, \s_\a)}
\def\qb{(A_\b, G_\b, \s_\b)}
\def\Fin{\mathop{\rm Fin}\nolimits}
\def\br{\mathop{\rm br}\nolimits}
\def\fin{\mathop{\rm fin}\nolimits}
\def\Ann{\mathop{\rm Ann}\nolimits}
\def\mspec{\mathop{\rm mspec}\nolimits}
\def\spec{\mathop{\rm spec}\nolimits}
\def\End{\mathop{\rm End}\nolimits}
\def\Mod{\mathop{\rm Mod}\nolimits}
\def\bfb{\mathop{\rm\bf b}\nolimits}
\def\bfd{\mathop{\rm\bf d}\nolimits}
\def\bfi{\mathop{\rm\bf i}\nolimits}
\def\bfj{\mathop{\rm\bf j}\nolimits}
\def\tr{{\rm tr}}
\def\df{{\rm df}}
\def\bfk{\mathop{\rm\bf k}\nolimits}
\def\bfc{\mathop{\rm\bf c}\nolimits}
\def\bEnd{\mathop{\rm\bf End}\nolimits}
\def\id{\mathop{\rm id}\nolimits}
\def\Ass{\mathop{\rm Ass}\nolimits}
\def\Ext{\mathop{\rm Ext}\nolimits}
\def\Ines{\mathop{\rm Ines}\nolimits}
\def\Hom{\mathop{\rm Hom}\nolimits}
\def\Aut{\mathop{\rm Aut}\nolimits}
\def\iso{\mathop{\rm Iso}\nolimits}
\def\bHom{\mathop{\rm\bf Hom}\nolimits}
\def\Rk{ R_\k-\mathop{\bf Mod}}
\def\Rn{ R_n-\mathop{\bf Mod}}
\def\map{\mathop{\rm map}\nolimits}
\def\cf{\mathop{\rm cf}\nolimits}
\def\msp{\mathop{\rm msp}\nolimits}
\def\top{\mathop{\rm top}\nolimits}
\def\product{\mathop{\rm prod}\nolimits}
\def\ker{\mathop{\rm Ker}\nolimits}
\def\Bext{\mathop{\rm Bext}\nolimits}
\def\Br{\mathop{\rm Br}\nolimits}
\def\dom{\mathop{\rm Dom}\nolimits}
\def\min{\mathop{\rm min}\nolimits}
\def\im{\mathop{\rm Im}\nolimits}
\def\max{\mathop{\rm max}\nolimits}
\def\AX{A[X]}
\def\rk{\mathop{\rm rk}}
\def\Diam{\diamondsuit}
\def\E{{\mathbb E}}
\def\Z{{\mathbb Z}}
\def\Q{{\mathbb Q}}
\def\N{{\mathbb N}}
\def\bQ{{\bf Q}}
\def\bF{{\bf F}}
\def\bX{{\bf X}}
\def\bS{{\mathbb S}}
\def\AA{{\cal A}}
\def\BB{{\cal B}}
\def\CC{{\cal C}}
\def\DD{{\cal D}}
\def\TT{{\cal T}}
\def\FF{{\cal F}}
\def\GG{{\cal G}}
\def\PP{{\cal P}}
\def\SS{{\cal S}}
\def\XX{{\cal X}}
\def\YY{{\cal Y}}
\def\fAb{\mathfrak{ Ab}}
\def\fA{\mathfrak{ A}}
\def\fS{{\mathfrak S}}
\def\fH{{\mathfrak H}}
\def\fz{{\mathfrak z}}
\def\fU{{\mathfrak U}}
\def\fV{{\mathfrak V}}
\def\fW{{\mathfrak W}}
\def\fK{{\mathfrak K}}
\def\fp{{\mathfrak p}}
\def\PT{{\mathfrak{PT}}}
\def\T{{\mathfrak{T}}}
\def\fX{{\mathfrak X}}
\def\fP{{\mathfrak P}}
\def\X{{\mathfrak X}}
\def\Y{{\mathfrak Y}}
\def\F{{\mathfrak F}}
\def\C{{\mathfrak C}}
\def\B{{\mathfrak B}}
\def\J{{\mathfrak J}}
\def\fN{{\mathfrak N}}
\def\fM{{\mathfrak M}}
\def\Fk{{\F_\k}}
\def\bar{\overline }
\def\Bbar{\bar B}
\def\Cbar{\bar C}
\def\Pbar{\bar P}
\def\Tbar{\bar T}
\def\Dbar{\bar D}
\def\abar{\bar a}
\def\bbar{\bar b}
\def\fbar{\bar f}
\def\ybar{\bar y}
\def\vbar{\bar \v}
\def\nubar{\bar \nu}
\def\Abar{\bar A}
\def\a{\alpha}
\def\e{\varepsilon}
\def\f{\phi}
\def\o{\omega}
\def\k{\kappa}
\def\hs{\widehat\s}
\def\hv{\widehat\v}
\def\hF{\widehat F}
\def\v{\varphi}
\def\z{\zeta}
\def\s{\sigma}
\def\l{\lambda}
\def\m{\mu}
\def\lo{\l^{\aln}}
\def\b{\beta}
\def\g{\gamma}
\def\d{\delta}
\def\al{\aleph}
\def\ale{\aleph_1}
\def\aln{\aleph_0}
\def\alm{\aleph_{\it{fms}}}
\def\alfwi{\aleph_{\it{fwi}}}
\def\alfsi{\aleph_{\it{fsi}}}
\def\alwc{\aleph_{\it{fwc}}}
\def\Cont{2^{\aln}}
\def\nld{{}^{ n \downarrow }\l}
\def\n+1d{{}^{ n+1 \downarrow }\l}
\def\hsupp#1{[[\,#1\,]]}
\def\size#1{\left|\,#1\,\right|}
\def\Binfhat{\widehat {B_{\infty}}}
\def\Z{{\mathbb Z}}
\def\Q{{\mathbb Q}}
\def\Zhat{\widehat \Z}
\def\Mhat{\widehat M}
\def\Rhat{\widehat R}
\def\Phat{\widehat P}
\def\Fhat{\widehat F}
\def\fhat{\widehat f}
\def\Ahat{\widehat A}
\def\Chat{\widehat C}
\def\Ghat{\widehat G}
\def\Bhat{\widehat B}
\def\shat{\widehat s}
\def\Btilde{\widetilde B}
\def\Ftilde{\widetilde F}
\def\restr{\mathop{\upharpoonright}}
\def\to{\rightarrow}
\def\arr{\longrightarrow}
\def\TLC{$\Z$-{\rm TLC}}
\newcommand{\norm}[1]{\text{$\parallel\! #1 \!\parallel$}}
\newcommand{\supp}[1]{\text{$\left[ \, #1\, \right]$}}
\def\set#1{\left\{\,#1\,\right\}}
\newcommand{\mb}{\mathbf}
\newcommand{\cp}{\widehat}
\newcommand{\wt}{\widetilde}
\newcommand{\dsum}{\bigoplus}
\newcommand{\norma}[1]{\mbox{$\parallel\! #1 \!\parallel_A$}}
\newcommand{\suppx}[1]{\mbox{$\left[ \, #1\, \right]_X$}}
\newcommand{\suppa}[1]{\mbox{$\left[ \, #1\, \right]_A$}}
\newcommand{\suppl}[1]{\mbox{$\left[ \, #1\, \right]_\l$}}
\newcommand{\card}[1]{\mbox{$\left| #1 \right|$}}
\newcommand{\union}{\bigcup}
\newcommand{\inters}{\bigcap}
\newcommand{\pure}{\subseteq_\ast}
\newcommand{\mapr}[1]{\xrightarrow{#1}}

\renewcommand{\labelenumi}{(\roman{enumi})}
\newcommand{\dach}[1]{\hat{\vphantom{#1}}}
\numberwithin{equation}{section}
\def\Pf{\smallskip\goodbreak{\sl Proof. }}
\def\Fp{\hfill \lower.3333ex\hbox{\vbox{\hrule\hbox{\vrule\phantom{\vrule
height 6.83333pt depth 1.94444pt width 8.77777pt}\vrule}\hrule}}
\ifmmode\let\next\relax\else\let\next\par\fi \next}
\def\less{\smallsetminus}
\def\Set#1{\{#1\}}
\def\Tree{{}^{\rho>{}}2}
\def\To#1{\buildrel #1 \over \longrightarrow}
\def\IM{\mathop{\text{Im}}}
\def\norm#1{\|#1\|}
\catcode`\@=11
\def\bigcupdot{\dot{\DOTSB\bigcup@\slimits@}}
\def\rprod_#1{{\prod_{#1}\nolimits^{\mathbf{red}}}}
\def\iff{\mathrel{\Longleftrightarrow}}

\title{On Crawley Modules
\thanks{This work is supported by the project No.
I-706-54.6/2001 of the German-Israeli Foundation for Scientific
Research \& Development\newline The authors would also like to
thank the Edmund Landau Center for partial support of this
research \newline Shelah's list of publications GbSh 833.
\newline\indent
AMS-subject classification 20K40, 03E50, 13C10.}
\author{R\"udiger G\"obel \\FB 6, Mathematik,
Universit\"at Duisburg Essen, \\ 45117 Essen, Germany
\\e-mail: \ \ {\small R.Goebel@Uni-Essen.De}\\ and\\
Saharon Shelah\\
Institute of Mathematics,
Hebrew University, Jerusalem, Israel \\
and Rutgers University, New Brunswick, NJ, U.S.A \\
{\small e-mail: shelah@math.huji.ac.il}}}
\date{} \maketitle

\begin{abstract}
This continues recent work in a paper by Corner, G\"obel and
Goldsmith \cite{CGG}. A particular question was left open: Is it
possible to carry over the results concerning the undecidability
of torsion--free Crawley groups to modules over the ring of
$p$-adic integers? We will confirm this and also strengthen one of
the results in \cite{CGG} by replacing the hypothesis of
$\Diamond$ by $CH$. For details see the introduction.
\end{abstract}

\section{Introduction}

This note is an extension of the recent work in \cite{CGG}. Let
$R$ be a principal ideal domain with quotient field $Q$. In this
paper we will consider torsion--free $R$--modules. Crawley
$p$-groups were investigated intensively in the last two decades
of the last century. Megibben \cite{Me} showed that parallel to
Whitehead groups, the existence of `proper' Crawley groups (those
which are not direct sums of cyclics) is undecidable. Further
fundamental results for Crawley $p$-groups were derived in two
adjoint papers by  Mekler, Shelah \cite{MS1,MS2}. Parallel to the
case of abelian Crawley $p$--groups we can define Crawley modules
in the torsion--free case as outlined in \cite{CGG}.

\begin{definition} An $R$-module $G$ is a Crawley module if for
any pair $M,N$ of pure and dense submodules of corank $1$ there is
an automorphism $\Theta\in \Aut_RG$ with $M\Theta =
N$.\end{definition}

The hypothesis on $M$ is equivalent to say that $G/M\cong Q$.

An old theorem of Jeno Erd\"os (see the early edition of Fuchs
\cite{Fu0} and \cite{CGG}) can be reformulated to give a first
result on torsion--free Crawley modules.

\begin{theorem} \label{erdo} Free $R$-modules are Crawley modules.\end{theorem}

Moreover, the problem whether $\ale$-free Crawley groups ($R =
\Z$) are free is undecidable. This follows from the main results
of \cite{CGG}:

\begin{theorem} \label{cgg} Let $R$ be a principal ideal domain and
$\k$ a regular cardinal $ >\size{R}$.
\begin{enumerate}
    \item {\rm(Assume $\Diamond_\k(E)$ for all stationary sets
    $E\subseteq \k$.)} Any $\k$--free $R$--Crawley module of
    cardinality $\k$ is free.
    \item {\rm (Assume ZFC + MA.)} If $\size{R}<\Cont$, then any
    $\ale$-separable, $\ale$-coseparable $R$-module of rank $\ale$
    is a Crawley module. Thus there are non-free Crawley
    $R$-modules of rank $\ale$ if $R$ is not a field and
    $\ale < \Cont$.
\end{enumerate}
\end{theorem}

We just note, that any strongly $\ale$-free $R$-module of
cardinality $\ale$ satisfies the hypothesis in Theorem \ref{cgg}
$(ii)$, thus is a Crawley module. The existence of non-free but
strongly $\ale$-free $R$-modules of cardinality $\ale$ under (ZFC
+ MA + $\ale < \Cont$) is well known, see \cite{EM} or directly
Shelah \cite{Sh0}. We want to replace $\Diamond_\k(E)$ for $\k
=\ale$ in Theorem \ref{cgg} by $\Cont < 2^{\ale}$ and add
$\size{R}\le \Cont$ in $(i)$ and in $(ii)$, respectively; see the
results in the appropriate sections. The idea of the proof of the
second assertion may also be useful for other applications. In
particular we will derive the following new

\begin{corollary} \label{intrcor} Let $R= J_p$ be the ring of $p$-adic integers.
\begin{enumerate}
    \item If $\Cont < 2^{\ale}$, then any reduced, torsion-free
    Crawley $R$-module of rank at most $\ale$ is free.
    \item {\rm (ZFC + MA + $\ale < \Cont$)} There are non-free,
    torsion-free, reduced Crawley $R$-modules of rank $\ale$.
\end{enumerate}
\end{corollary}
\bigskip

\section{The Construction using the Weak Diamond.}

The notion of the weak diamond principle $(\Phi)$ which is
equivalent to $2^{\aln} < 2^{\ale}$ (a weak form of the special
continuum hypothesis) comes from \cite{DS}; it is stated in Eklof,
Mekler \cite[pp. 147 -- 152]{EM}.

$\Phi_\l(S):$ \ {\it Let $S$ be a stationary subset of an
uncountable, regular cardinal $\l$ (thus $\cf \l =\l$) and let $D=
\{ f\in {}^\a 2: \a <\l\}$ (where $2=\{0,1\}$). Then $S$ has the
weak diamond property $\Phi_\l(S)$ if the following holds.

For any (coloring) function $c: D\arr 2 =\{0,1\}$, there is a weak
diamond function $\eta\in {}^\l 2$ (for $c$) such that for all
$f:\l\arr 2$ the set $\{\d\in S: (f\restr \d)c = \d\eta\}$ is
stationary in $\l$.}
\bigskip

Here we will use a stronger version which still is equivalent to
the weak diamond property.
\bigskip

\begin{lemma}\label{diam} Let $S$ be a stationary subset of an
uncountable, regular cardinal $\l$. Then $S$ has the weak diamond
property $\Phi_\l(S)$ if and only if the following holds for the
cardinal $\mu =: 2^{<\l}$ and the set $D= \{ f\in {}^\a\mu: \a <\l\}$.\\
If $c: D\arr 2$ is a coloring of $D$, then for every function
$f:\l\arr \mu$ there is a weak diamond sequence $\eta\in {}^\l 2$
(depending on $c$) such that $\{\d\in S: (f\restr \d)c = \d\eta\}$
is stationary in $\l$.
\end{lemma}
\bigskip

\Pf This is a special case of Shelah \cite[Appendix, Theorem 1.10,
see Section 3]{Sh:f}. \Fp

\bigskip
We have a particular case for $\l =\ale$.

\bigskip
If $S$ is a stationary subset of $\ale$, $\mu =: 2^{\aln} <
2^{\ale}$ and $c:\{f\in {}^\a\mu: \a < \o_1\}\arr 2$, then there
is $\eta:\o_1\arr 2$ such that for all $f:\o_1\arr \mu$ the set
$\{\d\in S: (f\restr \d)c = \d\eta\}\subseteq \o_1$ is stationary.

\bigskip

We derive from $\Phi_\l(S)$ a lemma which can be applied
immediately for Crawley modules.

\begin{lemma}\label{combi} Let $S$ be a stationary subset of a regular,
uncountable cardinal $\l$ such that $\Phi_\l(S)$ holds. Suppose
that $A$ is a set of cardinality $\card{A}=2^{<\l}$ and for any
$\eta\in {}^\l 2$ there are a cub $C_\eta$ of $\l$ and a function
$\Theta_\eta: \l \arr A$, then we can find $\a\in\l$ and
$\eta_0,\eta_1 \in {}^\l 2 \ (i=0,1)$ such that the following
holds.
\begin{enumerate}
    \item $\eta_0\restr \a =\eta_1\restr \a$
    \item $\a\eta_0=0$ and $\a\eta_1=1$
    \item $\Theta_{\eta_0}\restr \a=\Theta_{\eta_1}\restr \a$
    \item $\a\in C_{\eta_0}$
\end{enumerate}
\end{lemma}

\Pf Let $\mu = 2^{<\l}$ and choose a one to one (coding) function
$g:A\times \{0,1\}\arr \mu$. This is used to define for each
$\eta\in {}^\l 2$ a function $f_\eta: \l\arr \mu$ by
 $$(i \Theta_\eta, i\eta)g =i f_\eta \text{ for every } i\in \l.$$

In order to apply $\Phi_\l(S)$ we also define a coloring $c: D
\arr 2$, where $D= \{f:\a\arr \mu, \a\in\l\}$. If $f:\a\arr \mu$
is in $D$, then define

\begin{align} \label{color}
    f c  = \left\{
    \begin{array}{c@{\qquad}l}
    1 &\text{ if there is  } \eta\in {}^\l 2 \text{ such that }
    f_\eta\restr\a= f,  \a \in C_\eta \text{ and } \a\eta=0\\
    0 &\text{otherwise.}
    \end{array} \right.
\end{align}

If $\eta$ is as in the first line of (\ref{color}) we call $\eta$
{\em a witness for } $f$. Now let $\eta\in {}^\l 2$ be a weak
diamond function for $c$ given by $\Phi_\l(S)$. Hence $S_1 =
\{\d\in S: (f_\eta\restr \d)c = \d \eta\}$ is a stationary subset
of $\l$ and since $C_\eta$ is a cub, also $S_2 = S_1\cap C_\eta$
is stationary in $\l$. Hence clearly $S_2\ne \emptyset$ and there
is some $\d\in S_2$.

We first claim that

\begin{eqnarray} \text{ If } \d\in S_2 \text{ and } \eta \text{ is
a weak diamond function, then } \d\eta = 1.
\end{eqnarray}

\relax From $ \d\in S_2$ follows $(f_\eta\restr \d)c = \d\eta$ and $\d\in
C_\eta$, moreover $f_\eta$ extends $f_\eta\restr\d$. If also
$\d\eta =0$, then $(f_\eta\restr \d)c =0$ by the last observation
but also $(f_\eta\restr \d)c =1$ by (\ref{color}), which is a
contradiction and necessarily $\d\eta =1$.

We now choose $\eta_1=\eta$ for such a $\d$ and by the last claim
follows $\d\eta=1$. Now the first line of (\ref{color}) gives us a
witness $\eta_0$ for $f_{\eta_1}\restr \d$. By the first line of
(\ref{color}) follows $\d\eta_0=0$ and $\d\in C_{\eta_0}$. Hence
$(ii)$ and $(iv)$ are shown.

Then $  f_{\eta_0}\restr \d = f_{\eta_1}\restr \d$, hence $(i
\Theta_{\eta_0}, i\eta_0)g = (i \Theta_{\eta_1}, i\eta_1)g$ for
all $i<\d$. We apply $g^{-1}$ to get $\Theta_{\eta_0}\restr \d =
\Theta_{\eta_1}\restr \d$ and $\eta_0\restr \d =\eta_1\restr \d$.
Thus (i) and (iii) follow as well. \Fp
\bigskip

We now return to Crawley modules. If $\k$ is a regular cardinal,
then any $\k$-free $R$-module of rank $\k$ has a well defined
$\Gamma$-invariant which is an equivalence class of a subset $E$
of $\k$ (i.e. $E\equiv E' \iff E\cap C = E' \cap C$ for some cub
$C$ in $\k$), see Eklof, Mekler \cite[p. 86]{EM}.

Here we want to show the following

\begin{theorem}\label{diamthm} Let $R$ be a principal ideal domain and let $G$ be
an $\ale$-free $R$-module of rank $\ale$ with $\Gamma$-invariant
$\Gamma(E)$ coming from a stationary set $E\subseteq \ale$. If the
weak diamond holds for $E$, then $G$ is not a Crawley module.\\
Thus, assuming $\Cont < 2^{\ale}$, an $\ale$-free $R$-module of
rank $\ale$ is Crawley if and only if it is free.\end{theorem}

\Pf Let $G= \bigcup\limits_{\a\in\ale} G_\a$ be an
$\ale$-filtration of $G$ with $G_\a$ pure in $G$ of countable rank
for each $\a\in\o_1$. Then $$E =\{\a\in \o_1, G/G_\a \text{ is not
} \ale-\text{free}\}.$$

If $E$ is not stationary, then $G$ is free and the claim follows
from Erd\"os's Theorem \ref{erdo}. We may assume for contradiction
that $E$ is stationary and a set of limit ordinals without
restrictions. If $\a\in E$, then $G_{\a+1}/G_\a$ is countable and
not free. By Stein's theorem (see Fuchs \cite{Fu}) we can
decompose $G_{\a+1}/G_\a$ into a direct sum of a free module and a
complement different from $0$ with trivial dual (having no proper
homomorphisms into $R$). We can absorb the free summand into the
next $G_{\a+2}$ and thus assume that it is $0$. Hence
$G_{\a+1}/G_\a\ne 0$ and  $(G_{\a+1}/G_\a)^* =0$ ($G_{\a+1}/G_\a$
is not free of minimal rank $n_\a\le \o$). There is a free
$R$-module $X'\subseteq G_{\a+1}/G_\a$ such that
$X'_*=G_{\a+1}/G_\a$. Its preimage $X\subseteq G_{\a+1}$ is of the
form
\begin{eqnarray}\label{nonfree} X = \bigoplus\limits_{m < n_\a} Rx_{\a m}\subseteq
G_{\a+1} \text{ with } G_{\a+1} = (G_\a \oplus X)_* \text{ and }
(G_{\a+1}/G_\a)^* = (X'_*)^* =0.
\end{eqnarray}

We choose one of these pure and dense maximal submodules
$M\subseteq G$ of corank $1$ and create a large family of its
relatives:

By definition of $M$ we can find $z\in G\setminus M$ such that $G
= \langle M, Rz\rangle_*$. Moreover (changing the filtration) we
may assume that $z\in G_0$. By induction on $\a\in\o_1$ we
determine for each $\eta\in {}^\a2$ an epimorphism $$\v_\eta:
G_\a\arr Q$$.

If $\a =0$, then $\eta =\langle \rangle$ is the empty map, we send
$z$ to $1\in Q$. We also assume that $G_0$ has infinite rank and
map infinitely many free generators of $G_0$ onto $\Q$. By
injectivity this map extends further to an epimorphism
$\v_{\langle \rangle}: G_0 \arr Q$. If $\nu \trianglelefteq \eta$
is an initial segment, then we want $\v_\nu \subseteq \v_\eta$.
Suppose that $\v_\eta$ is constructed for all $\eta\in {}^{\a>}2$
We distinguish three cases. If $\a$ is a limit ordinal, and
$\eta\in {}^\a2$, then let $\v_\eta = \bigcup\limits_{\b<\a}
\v_{\eta\restr \b}$. If $\b\in \o_1\setminus E$ and $\a = \b +1$ ,
then we assume that $\eta$ takes only one value at $\a$, say
$\a\eta = 0$. Now any $\eta \in {}^\b2$ has a unique extension
$\nu =\eta ^\wedge \langle 0\rangle \in {}^\a2$. Since $G_\b$ is a
summand of $G_{\a}$ with free complement, we can extend $\v_\eta$
arbitrarily to $\v_\nu: G_\a\arr Q$. Finally we assume that $\b\in
E$ and select the elements $x_{\b m} \ (m\le n_\b)$ by the choice
of $X$ in (\ref{nonfree}). If $\eta \in {}^\b2$ and $\v_\eta:
G_\b\arr Q$ is given, then we want to choose $\v_{\eta ^\wedge
\langle 0\rangle}$ and $\v_{\eta ^\wedge \langle 1\rangle}$ such
that

\begin{eqnarray}\label{equ}
 x_{\b m}\v_{\eta ^\wedge \langle 1\rangle} - x_{\b m}\v_{\eta
^\wedge \langle 0\rangle}= 1 \text{ for all } m\le
n_\b.\end{eqnarray}

Using (\ref{nonfree}) write $G_\a =G_\b\oplus X$ for $X =
\bigoplus\limits_{m\le n_\b} Rx_{\b m}$ and we can extend
$\v_\eta$ to both of $\v'_\eta, \v''_\eta\in \Hom(G_\b\oplus X,Q)$
with $x_{\b m}\v'_\eta = 0$ and $x_{\b m}\v''_\eta = 1$ for all
$m\le n_\b$. These maps can be extended further to $\v_{\eta
^\wedge \langle 0\rangle}, \v_{\eta ^\wedge \langle 1\rangle}\in
\Hom(G_\a,Q)$ by injectivity. We have $x_{\b m}\v_{\eta ^\wedge
\langle 0\rangle}=0$ and $x_{\b m}\v_{\eta ^\wedge \langle
1\rangle}=1$ for all $m\le n_\b$ and in particular (\ref{equ})
holds.

This finishes the inductive construction. If $\eta \in {}^{\o_1}2$
(with $\eta \restr \o_1\setminus E$ the zero map), then we obtain
the relatives of $M$:
$$\v_\eta = \bigcup\limits_{\a\in\o_1}\v_{\eta\restr \a} \text{
and } N_\eta = \ker\v_\eta.$$

Hence $G/\ker \v_\eta\cong Q$ and $\{ N_\eta : \eta\in
{}^{\o_1}2\}$ is a family of pure and dense submodules of corank
$1$.

We claim that there is $\eta\in {}^{\o_1}2$ such that no
$\Theta\in \Aut_RG$ satisfies $M\Theta = N_\eta$. From this the
Theorem follows.

Assume for contradiction that for any $\eta\in {}^{\o_1}2$ there
is $\Theta_\eta\in \Aut_RG$ with $M\Theta_\eta = N_\eta$. Now we
specify the choice of the cubs $C_{\eta}\subseteq \l$ for Lemma
\ref{combi}: By the usual back and forth argument we may assume
that $G_\a\Theta_{\eta}=G_\a$ for all $\a\in C_{\eta}$. By the
weak diamond and Lemma \ref{combi} there are some $\eta_0,\eta_1
\in {}^{\o_1}2$ and $\a\in\o_1$ such that
\begin{enumerate}
    \item $\a$ is the branch point of $\eta_0$ and $\eta_1$:
    $\eta_0\restr \a = \eta_1\restr \a =\nu$ but $\a \eta_0 \ne \a\eta_1$
    \item  $\a\eta_0= 0$ and $\a\eta_1= 1$,
    \item  $\Theta_{\eta_0}\restr G_\a = \Theta_{\eta_1}\restr G_\a$.
    \item  $\a \in C_{\eta_0}.$
\end{enumerate}
If we put $\Theta =\Theta_{\eta_1}-\Theta_{\eta_0}$, then
$\Theta\in \End G$ induces $\theta: G/G_\a\arr G$. But necessarily
$\a\in E$, hence $(G_{\a+1}/G_\a)^*=0$ by (\ref{nonfree}). On the
other hand $G$ is $\ale$-free and $(G_{\a+1}/G_\a)$ is countable,
hence $G_{\a+1}\Theta = 0$. It follows that
\begin{eqnarray}\label{fix} \Theta_{\eta_1}\restr G_{\a+1} =
\Theta_{\eta_0}\restr G_{\a+1}.\end{eqnarray}

Since $\a\eta_0=0$ by (ii), we have $\eta_0=\nu ^\wedge \langle
0\rangle$ and $x_{\a 0} \v_{\nu^\wedge \langle 0\rangle} = 0$ by
definition of the extension map. But $\v_{\nu^\wedge \langle
0\rangle}\subseteq \v_{\eta_0}$ and so $x_{\a 0} \v_{\eta_0}=0$.
We derive $x_{\a 0}\in N_{\eta_0}$.

By the choice of $C_{\eta_0}$ it will follow
 $y^*\Theta_{\eta_0}\in G_{\a+1}$: Note that $\Theta_{\eta_0}^{-1}$
leaves $G_\a$ invariant and thus also maps $G_{\a+1}$ into itself
because $G/G_{\a+1}$ is $\ale$-free but $(G_{\a+1}/G_\a)^*= 0$ by
(\ref{nonfree}). (We say that $G_{\a+1}$ is the Chase radical of
$G$ over $G_\a$.) From $y^*\Theta_{\eta_0}\in G_{\a+1}$ and
(\ref{fix}) also follows $x_{\a 0} = y^*\Theta_{\eta_0}$. But now
$x_{\a 0}= y^*\Theta_{\eta_1}\in G_{\a+1}\cap N_{\eta_1} = \ker
\v_{\nu ^\wedge\langle 1\rangle}$, so $x_{\a 0}\v_{\nu
^\wedge\langle 1\rangle}= 0$. Hence $x_{\a 0}\v_{\nu
^\wedge\langle 0\rangle}- x_{\a 0}\v_{\nu ^\wedge\langle
1\rangle}=0$, which contradicts (\ref{equ}). \Fp

\bigskip
The proof immediately gives the a stronger result:

If $M\subseteq G$ is a pure and dense submodule of corank $1$, the
we denote by $[M]= M\Aut_RG$ the orbit of $M$ under the action of
the automorphism group.

\begin{corollary} {\rm (ZFC + $2^{\aln}< 2^{\ale}$):}
If $G$ is a module as in the Theorem \ref{diamthm} which is not
free and of cardinality $\ale$, then $G$ has $\ale$ distinct
orbits of pure, dense submodules of corank $1$.\end{corollary}

\Pf Partition $E$ into $\ale$ stationary sets. \Fp

\section{Crawley modules under Martin's Axiom}\label{models}

The proof of our next theorem is inspired by the
L\"owenheim-Skolem argument. We first state and prove an
observation mentioned by Brendan Goldsmith. This will be used to
reduce the theorem to a result in \cite{CGG}.

\begin{observation} \label{obs1} Every subring of a principal
ideal domain $R$ is a subring of principal ideal domain contained
in $R$ and of the same cardinality.
\end{observation}
\Pf Let $R_1\subseteq R$ be the given rings, let $Q$ be the field
of fractions of $R_1$ and set $R_2= R \cap Q$ which is a subring
of $R$ of cardinality $|R_1|$ containing $R_1$. We will show that
$R_2$ is a principal ideal domain. If $I$ is an ideal of $R_2$,
then $IR$ is an ideal of $R$ and so $IR =aR$ for some $a \in R$.
However $a \in IR$ and so $a = ir$ for some $i \in I, r \in R$.
Now if $x$ is an arbitrary element of $I$, then $ x = ar_x$ for
some $r_x \in R$ and so $x= i r r_x = i t_x$ say where $t_x = r
r_x \in R$. However $t_x = x/i \in K \cap R = R_2$. So $x = i t_x$
is a product of an element in $I$ and in $R_2$. Since $x$ was
arbitrary in $I$ we have $I \subseteq i R_2$ and since the reverse
inclusion is trivial, we deduce that $I = i R_2$ is principal. \Fp

\begin{theorem} \label{mathm} {\rm (ZFC + MA)} Let $R$ be a principal
ideal domain. Any strongly $\ale$-free $R$-module of rank $\ale$
is a Crawley module.
\end{theorem}

\Pf Let $N_1,N_2\subseteq G$ be two pure and dense submodules of
corank $1$ of the strongly $\ale$-free $R$-module $G$ of rank
$\ale$. Moreover let $\v_i:G\arr Q$ be the canonical epimorphisms
with kernel $\ker \v_i = N_i \ (i=1,2)$. Choose an
$\ale$-filtration $G=\bigcup\limits_{\a\in\o_1} G_\a$ and let
$N_{i\a} = N_i\cap G_\a$ be the induced $\ale$-filtration on $N_i
\ (i=1,2)$. We can assume that there are $z_i\in G_0$ with
$z_i\v_i = 1\in Q$ for $i=1,2$. Moreover, let each $G_0,N_{i 0}$
have infinite rank (w.l.o.g.). Also choose an $R$-basis for each
$G_\a, N_{i\a}$, thus
$$ G_\a =\bigoplus\limits_{n\in\o} Rx^0_{\a n},\ N_{i\a} =\bigoplus
\limits_{n\in\o} Rx^i_{\a n} \text{ for } i= 1,2.$$

If $\a < \b\in\o_1, n\in\o$, then $$x^i_{\a n} =
\sum\limits_{m\in\o} a^{i\a}_{\b nm} x^i_{\b m} \text{ for some }
a^{i\a}_{\b nm}\in R.$$

Since $z_i\in G_0$, also $z_i= \sum\limits_{m\in\o}b^i_m x^0_{0
m}$ for some $b^i_m \in R$. Now we can choose a subring $R'$ of
$R$ of cardinality $\ale$ containing all these coefficients
$a^{i\a}_{\b nm}, b^i_m$.

Let $$G'_\a=\bigoplus\limits_{n\in\o} R'x^0_{\a n},\ N'_{i\a}
=\bigoplus \limits_{n\in\o} R'x^i_{\a n} \text{ for } i= 1,2$$

be the corresponding free $R'$-module. We have $\ale$-filtrations
$G'=\bigcup\limits_{\a\in\o_1} G'_\a$ and $N'_i
=\bigcup\limits_{\a\in\o_1}N'_{i\a}$; this uses that $R'$ is
sufficiently saturated! Moreover, passing from $G$ to $G'$ it is
easy to see that also $G'$ is a strongly $\ale$-free $R'$-module.
Also $z_i \in G'_0$ and $\v'_i = \v_i\restr G'$ is an
$R'$-homomorphism $G'\arr Q$.

By another back and forth argument, using that $\v_i: G\arr Q$ are
epimorphisms (so enlarging $R'$) we may assume that $\v_i'$ maps
$G'$ onto the quotient field $Q'$ of $R'$. By Observation
\ref{obs1} we also assume that $R'$ is a principal ideal domain.
Thus $G'/N'_i\cong Q'$ for $i=1,2$, and we can apply Theorem
\ref{cgg}$(ii)$. The $R'$-module $G'$ is a Crawley $R'$-module and
there is $\Theta'\in \Aut_{R'}G'$ with $N'_1\Theta'= N_2'$.

Finally we extend $\Theta'$ (uniquely) to $\Theta\in \Aut_RG$ such
that $N_1\Theta = N_2$:

There is a cub $C\subseteq \o_1$ such that $\Theta'_\a
=\Theta'\restr G'_\a\in \Aut_{R'}G'_\a$ for all $\a\in C$;
moreover $\Theta' = \bigcup\limits_{\a\in C} \Theta'_\a$ and $G' =
\bigcup\limits_{\a\in C} G'_\a$. These $R'$-automorphisms
$\Theta'_\a$ act on the free $R'$-module $G'_\a =
\bigoplus\limits_{n\in\o} R'x^0_{\a n}$; they extend naturally to
$R$-automorphisms $\Theta_\a$ of the free $R$-module $G_\a
=\bigoplus\limits_{n\in\o} Rx^0_{\a n}$ by tensoring with $R$.
Clearly $\Theta = \bigcup\limits_{\a\in C} \Theta_\a$ and
$N_1\Theta = N_2$.

(This also shows that $G$ is $\ale$-presented, i.e. it can be
expressed as the quotient of two free $R$-modules of rank $\ale$.)
\Fp
\bigskip

If $R$ is a principal ideal domain with quotient field $Q$
countably generated over $R$, then there is a multiplicatively
closed subset $\bS$ of $R\setminus \{0\}$ such that $\bS^{-1}R
=Q$. We may assume that $1\in \bS$ and choose an enumeration
$\{s_i : i\in\o\}=\bS$ with $s_0=1$. Moreover $q_n= \prod_{i\le
n}s_i$. In particular $0 =\bigcap\limits_{i\in\o} Rq_i$, so $R$ is
an $\bS$-ring; see \cite{GT}. Examples are the principal ideal
domains of $p$-adic integers $J_p$ with $\bS =\{p^n: n\in\o\}$ and
counter examples are polynomial rings $\Z[X]$ over $\Z$ in  sets X
of uncountably many, commuting variables. If $R$ is such an
$\bS$-ring, then we can construct Griffith's strongly $\ale$-free
$R$-module of rank $\ale$: Let $F =\bigoplus\limits_{i\in \o_1}
Re_i$ be the free $R$-module of rank $\ale$ and $\Fhat$ its
$\bS$-adic completion; note that the $\bS$-topology on $F$ is
Hausdorff by the above. Let $E\subseteq \o_1$ be a stationary
subset of limit ordinals and choose for any limit ordinal
$\d\in\o_1$ a ladder, i.e. a strictly increasing sequence of
successor ordinals $\d_n \ (n\in\o)$ with limit $\d$. We consider
the elements $$ v_\d^k = \sum_{n\ge k} (q_n q_k^{-1}) e_{\d_n} \in
\Fhat$$ and let
$$G=\langle F, Rv_\d^k : \d\in E, k\in \o\rangle$$ which is a pure submodule
of $\Fhat \cap \prod\limits_{i\in \o_1} Re_i$. It is easy to
check, that $G$ is strongly $\ale$-free with $\Gamma$-invariant
$\Gamma(E)\ne 0$, hence not free. Thus the next corollary follows
from Theorem~ \ref{mathm}.

\begin{corollary} \label{MAcor} {\rm (ZFC + MA + $\ale <\Cont$)} If $R$ is a
principal ideal domain with quotient field $Q\ne R$ countably
generated over $R$, then there are $\ale$-free but not free
Crawley $R$-modules of rank $\ale$.\end{corollary}

Corollary \ref{intrcor} follow immediately from Corollary
\ref{MAcor} and Theorem \ref{mathm}.
\goodbreak

\end{document}